\documentclass{article}
\usepackage{amsmath,amsfonts,amsthm}
\newtheorem{theorem}{Theorem}
\newtheorem{lemma}[theorem]{Lemma}

\newcommand\pto{\overset{\mathrm{p}}{\to}}
\newcommand\vv{{\bf v}}
\newcommand\vw{{\bf w}}
\newcommand\vx{{\bf x}}
\newcommand\vy{{\bf y}}
\newcommand\ZZ{{\mathbb Z}}
\newcommand\RR{{\mathbb R}}
\newcommand\Leb{\mu}
\newcommand\bpm{\mu}%
\newcommand\kk{\kappa}
\newcommand\sss{{\mathcal S}}
\newcommand\dd{\,d}
\newcommand\norm[1]{\ensuremath{\|#1\|}}
\newcommand\good{{Y}}
\newcommand\CC{{\mathcal C}}
\newcommand\bb[1]{\bigl(#1\bigr)}

\newcommand{\refT}[1]{Theorem~\ref{#1}}

\newcommand{\refL}[1]{Lemma~\ref{#1}}

\renewcommand\Pr{{\mathop{\mathbb P{}}\nolimits}}

\newcommand\kernelsTtwo{3.1}
\newcommand\kernelsTtwob{3.6}
\newcommand\kernelsRpath{2.8}

\newcommand\webcite[1]{\hfil\penalty0\texttt{\def~{\~{}}#1}\hfill\hfill}

\begin{document}
\title{Spread-out percolation in $\RR^d$}
\date{16th October 2006}

\author{B\'ela Bollob\'as\thanks{Department of Mathematical Sciences,
University of Memphis, Memphis TN 38152, USA}
\thanks{Trinity College, Cambridge CB2 1TQ, UK}
\thanks{Research supported in part by NSF grant ITR 0225610}
\and Svante Janson%
\thanks{Department of Mathematics, Uppsala University,
 PO Box 480, SE-751 06 Uppsala, Sweden}
\and Oliver Riordan%
\thanks{Royal Society Research Fellow, Department of Pure Mathematics
and Mathematical Statistics, University of Cambridge, Wilberforce Road,
Cambridge CB3 0WB, UK}}
\maketitle

\begin{abstract}
Fix $d\ge 2$, and let $X$ be either $\ZZ^d$ or the points of a Poisson process in
$\RR^d$ of intensity 1.  Given parameters $r$ and $p$, join each pair
of points of $X$ within distance $r$ independently with probability
$p$. This is the simplest case of a `spread-out' percolation model
studied by Penrose~\cite{Penrose}, who showed that, as $r\to\infty$,
the average degree of the corresponding random graph at the
percolation threshold tends to $1$, i.e., the percolation threshold
and the threshold for criticality of the naturally associated
branching process approach one another. Here we show that this result
follows immediately from of a general result of \cite{kernels} on
inhomogeneous random graphs.
\end{abstract}

\section{Introduction and results}

The study of percolation and the study of the emergence of the giant
component in a random graph are closely related topics. In both cases,
one can phrase the key question as follows: `As $n\to \infty$, for
what parameters does a certain $n$-vertex random graph have with high
probability a component of order $\Theta(n)$?'
The key difference is that in percolation there is some global
geometric structure: for example, the graph might be a random subgraph
of $\ZZ^d$, or of a finite portion of $\ZZ^d$,
or it might be the graph formed by the points of
a Poisson process in $\RR^d$, joined if they are within a certain distance $r$.
In the classical theory of random graphs, the graph model, $G(n,p)$, is
totally structureless: each pair of vertices is joined independently
with probability $p$, with $p=c/n$, $c$ constant, being the
appropriate normalization for the very simplest results on the giant
component.
Of course, each instance of $G(n,p)$ has a rich structure,
but the model does not.
Many authors have studied inhomogeneous random graphs in which the model
does have some structure, but the behaviour of these random graphs
is still much closer to the behaviour of $G(n,p)$ than to percolation.
In particular,
it was shown in~\cite{kernels} that, for a very general inhomogeneous
model with independence between the edges, the threshold for the
emergence of a giant component is given by the point at which a
certain (multi-type) branching process becomes critical, generalizing
the classical result for $G(n,p)$. In contrast, in percolation there
are only a few models where the exact threshold can be determined:
unless one of a few very special things happens, it seems to be
impossible to give a simple formula for the critical point.

Penrose~\cite{Penrose} determined the asymptotic behaviour
of the critical parameters for a certain natural family
of two-parameter percolation models with global
geometric structure. A special case of this result,
as in the present
abstract but with $X=\ZZ^d$ and distance replaced by $\ell_\infty$-distance,
was proved independently by Bollob\'as and Kohayakawa~\cite{BK}.
Our aim in this paper is to show that Penrose's result is a simple consequence
of the main (and very general) results of~\cite{kernels}.

Writing $\Leb$ for the Lebesgue measure on $\RR^d$,
we say that a set $A\subset \RR^d$ is a {\em $\Leb$-continuity set}
if $A$ is $\Leb$-measurable and $\Leb(\partial A)=0$, where
$\partial A$ is the boundary of $A$. Note that we take the
$d$-dimensional measure of the boundary.

Throughout, the vertex set of our (infinite) random graph will be a
deterministic or random set $X\subset \RR^d$, of `density' one. The
dimension $d\ge 2$ will be fixed throughout.  The natural examples are
$\ZZ^d$ (perhaps with each point displaced by a small random amount),
any other (suitably scaled) lattice, or a Poisson process
of intensity $1$. 
Formally, we require two assumptions on the
distribution of the random discrete set $X$: the {\em density assumption} that if
$A\subset \RR^d$ is a bounded $\Leb$-continuity set,
then for any sequence $\vv_r\in \RR^d$ we have
\[
 \frac{1}{r^d}\bigl|X \cap (rA+\vv_r)\bigr| \pto \Leb(A)
\]
as $r\to\infty$,
where $\pto$ denotes convergence in probability,
$|\cdot|$ denotes the cardinality of a set, and $rA+\vv_r=\{r\vx+\vv_r:\vx\in
A\}$.  The second assumption we require is that well separated
regions are independent: there is a constant $D$ such that whenever $A, B\subset \RR^d$
are measurable sets at Euclidean distance at least $D$ from each other,
then the set-valued random variables
$X\cap A$ and $X\cap B$ are independent.

To state our form of Penrose's result, we consider a function
$\phi:\RR^d\to[0,\infty)$ satisfying the following assumptions: $\phi$
is {\em symmetric}, in that $\phi(\vx)=\phi(-\vx)$, bounded, continuous almost everywhere,
and satisfies $\int_{\RR^d} \phi(\vx) \dd\Leb(\vx)=1$. In addition, for convenience
we shall assume that $\phi$ is strictly positive in a neighbourhood of the origin.
To obtain the example described in the abstract, we choose for $\phi$ the function
that is $1$ on a ball of volume 1 centred at the origin, and $0$ otherwise.

Given a `scale' $r$ and a `degree parameter' $\lambda>0$, we form the random graph $G=G(X)=G_{r,\lambda}(X)$
with vertex set $X$ as follows: given $X$, for each pair $\vx,\vy\in X$
join $\vx$ and $\vy$ with probability
\[
 \min\left\{\lambda r^{-d}\phi\left(\frac{\vx-\vy}{r}\right) , 1\right\},
\]
independently of all other pairs.
Note that $\int_{\RR^d} r^{-d}\phi(\vx/r)\dd\Leb(\vx)=1$, so, at least in the Poisson case,
the average degree of a vertex of $G_{r,\lambda}(X)$ tends to $\lambda$ as $r\to\infty$
with $\lambda$ fixed.

\begin{theorem}\label{th_1}
Let $X$ and $\phi$ satisfy the conditions above with $d\ge 2$, and let $\lambda>1$ be fixed.
If $r$ is large enough, then with probability $1$ the graph $G_{r,\lambda}(X)$ has an infinite component.
\end{theorem}

The basic idea is to show that the neighbourhood of a vertex of $G(X)$ is `tree-like', and can be
approximated by a Galton-Watson branching process where each particle
has a Poisson number of children with mean $\lambda$. We shall approximate the local
structure of $G_{r,\lambda}(X)$
using the results of \cite{kernels}, and then deal with the global structure using
the concept of $k$-independent percolation.

\section{Proofs}
Throughout this section, $d\ge 2$, $D$, $\phi$, and the distribution of $X$ will be fixed.
Much of the time,
we shall rescale the model in the following natural way: let $X/r=\{\vx/r:\vx\in X\}$.
The (rescaled) graph $G'(X)=G_{r,\lambda}'(X)$ has vertex set $X/r$, and each
pair $\vx,\vy\in X/r$ is joined independently with probability $\min\{\lambda r^{-d}\phi(\vx-\vy),1\}$.

To prove Theorem~\ref{th_1} we shall need two results. The
first is a simple observation concerning locally dependent percolation.

A {\em bond percolation measure} on $\ZZ^d$ is a measure
on the set of assignments of a {\em state}, {\em open} or {\em closed},
to each edge of $\ZZ^d$, the graph with vertex set $\ZZ^d$ in which vertices
at Euclidean distance $1$ are adjacent.
Such a measure is {\em $k$-independent} if,
for every pair $S$, $T$ of sets of edges of $\ZZ^d$ at graph
distance at least $k$, the states of the edges in $S$ are independent of
the states of the edges in $T$. When $k=1$, the separation
condition is exactly that no edge of $S$ shares a vertex with an edge of $T$.

Measures of this type arise very naturally in static renormalization arguments, and
comparisons between $k$-independent measures and product measures (or arguments amounting to such comparisons)
have been considered by many people; see
Liggett, Schonmann and Stacey~\cite{LSS} and the references therein.

\begin{lemma}\label{l_d}
Let $d\ge 2$ and $k\ge 1$ be fixed. There is a $p_0=p_0(k)<1$
such that in any $k$-independent bond percolation measure 
on $\ZZ^d$ satisfying the
additional condition that each edge is open with probability at least $p_0$, with probability
1 there is an infinite path consisting of open edges.
\end{lemma}

This result is a special case of the very general main result of \cite{LSS}; in the form
above, it is essentially trivial. Note that without loss of generality we may take $d=2$,
as $\ZZ^d$ contains $\ZZ^2$ as a subgraph.
Here the value of $p_0$ is irrelevant, but in many contexts this value is very important.
For $k=1$, the best bound known is due to Balister, Bollob\'as and Walters~\cite{BBW}, who showed that
one can take $p_0(1)=0.8639$.

The second result we shall need is a special case of the main result
of \cite{kernels}; to state this we recall some definitions
from~\cite{kernels}.

Let $\sss$ be a separable metric space and $\bpm$ a Borel measure on
$\sss$ with $0<\bpm(\sss)<\infty$.
In this paper, $\sss$ will be either a cube in $\RR^d$ of side-length $L$,
or the union of two such cubes sharing a face, and $\bpm$ will be Lebesgue measure.

Let $\rho>0$ be a `density parameter' that will tend to infinity.
Here we shall fix $\sss$ and take $\rho=r^d$ with $r\to\infty$.
For each $\rho$, let $V_\rho$ be a deterministic or random finite subset
of $\sss$. The triple $(\sss,\bpm,(V_\rho))$ forms a {\em generalized vertex space}
if $\sss$ and $\bpm$ satisfy the conditions above, and
\begin{equation}\label{gvxs}
 \frac{1}{\rho} |V_\rho\cap A| \pto \bpm(A)
\end{equation}
as $\rho\to\infty$, for each $\bpm$-continuity set $A\subset \sss$.
Here, we shall take $V_\rho=(X/r)\cap \sss$, so the condition
above is simply that
\[
 r^{-d}|(X/r)\cap A| \pto \Leb(A),
\]
which follows from our density assumption on $X$.

A {\em kernel} on $(\sss,\bpm)$ is a symmetric, non-negative, Borel-measurable
function on $\sss\times\sss$; we shall consider kernels $\kk=\kk_{\phi,\sss,\lambda}$
given by $\kk(\vx,\vy)=\lambda\phi(\vx-\vy)$ for $\vx,\vy\in \sss$.
Note that as $\phi$ is continuous almost everywhere on $\RR^d$, the kernel $\kk$ is continuous
almost everywhere on $\sss\times\sss$.

A kernel $\kk$ is {\em irreducible} if, whenever $A$ is a measurable
subset of $\sss$ with $\kk=0$ a.e. on $A\times(\sss\setminus A)$, then
$\bpm(A)=0$ or $\bpm(\sss\setminus A)=0$.  The assumption that $\phi$
is strictly positive in a neighbourhood of the origin ensures that the
kernel $\kk=\kk_{\phi,\sss,\lambda}$ is irreducible.
(To see this, pick a $\delta>0$ such
that $\phi(\vx)>0$ on $B_{2\delta}(0)=\{\vx: \norm{\vx}< 2\delta\}$.
If $A\subset\sss$ with $0<\mu(A)<\mu(\sss)$, then, since $\sss$ is
connected, there is some $\vx\in \sss$ such that both $A$ and
$\sss\setminus A$ meet $B_\delta(\vx)$ in sets of positive measure.
Since $\kk$ is positive on $B_\delta(\vx)\times B_\delta(\vx)$ it
follows that $\kk$ is irreducible.)

Given a kernel $\kk$ on $(\sss,\bpm)$, let $T_\kk$ be the integral
operator on $(\sss,\bpm)$ with kernel $\kk$, defined by
\[
 (T_\kk f)(x) = \int_\sss \kk(x,y)f(y)\dd\bpm(y),
\]
for any (measurable) function $f:\sss\to \RR$ such that the integral is defined
(finite or $+\infty$) for a.e. $x$. For the bounded $\kk$ we consider,
$T_\kk f$ is defined for every $f\in L^2=L^2(\sss,\bpm)$, and
the operator $T_\kk$ maps $L^2$ into itself.
Let $\norm{T_\kk}$ be the operator norm of $T_\kk:L^2\to L^2$.

Theorem~\ref{th_k} below is a special case of parts of
the main results, Theorem~\kernelsTtwo\ (part (iii))
and Theorem~\kernelsTtwob, of \cite{kernels}. In Theorem~\ref{th_k},
$G_\rho(\kk)$ denotes the random graph with vertex set $V_\rho$ where,
given $V_\rho$, each pair $\{\vx$, $\vy\}$ of vertices is joined
independently with probability
$\min\{\kk(\vx,\vy)/\rho,1\}$. We write $C_i(G)$ for the number
of vertices in the $i$th largest component of a graph $G$.

\begin{theorem}\label{th_k}
Let $(\sss,\bpm,(V_\rho))$ be a generalized vertex space, and let $\kk$ be an irreducible,
bounded,
almost everywhere continuous kernel on $(\sss,\bpm)$.
If $\norm{T_\kk}>1$, then $C_1(G_\rho(\kk))/\rho\pto a$ as $\rho\to\infty$,
for some constant $a=a(\sss,\bpm,\kk)>0$,
while $C_2(G_\rho(\kk))/\rho\pto 0$.
\end{theorem}

For comparison with the statement in \cite{kernels}, note that the
additional condition there, that $\kk$ be `graphical' on
$(\sss,\bpm,(V_\rho))$, is not needed here. Indeed, taking $A=\sss$ in
\eqref{gvxs}, we have $\rho^{-1}|V_\rho|\pto \mu(\sss)<\infty$. In
particular, $\rho^{-1}|V_\rho|\le 2\mu(\sss)$ with probability
$1-o(1)$. Redefining $V_\rho$ to be empty if this inequality is not
satisfied, which changes $G_\rho(\kk)$ on a set of measure $o(1)$ and
hence does not affect the conclusion of Theorem~\ref{th_k}, we obtain
a new vertex space with $\rho^{-1}|V_\rho|$ bounded. But now
convergence in probability in \eqref{gvxs} implies convergence of all
moments. As noted in \cite{kernels} (Remarks \kernelsRpath\ and
8.2), under this condition any bounded, almost everywhere continuous
kernel is graphical.  [This argument, which applies to all `with
probability $1-o(1)$' results in \cite{kernels}, shows that the
definition of graphical there should perhaps be modified not to refer
to expectation; this is purely a matter of convenience, since one can
always modify the model on events with probability $o(1)$ as here.]

\begin{proof}[Proof of Theorem~\ref{th_1}]
Fix $\lambda>1$ throughout.
Let us rescale the vertex set $X$ as above, considering the graph $G'(X)$ with vertex set $X/r$.

Let $L$ be a (large) constant to be chosen later. Let $\sss_1=[0,L)^d$
and $\sss_2=[0,2L)\times [0,L)^{d-1}$. With $\sss=\sss_1$ or $\sss=\sss_2$,
let $\kk=\kk_{\phi,\sss,\lambda}$ be defined as above, and let
$f$ be the constant function on $\sss$ taking value $1$.
Then 
\begin{eqnarray*}
 (T_\kk f)(\vx) &=&
  \int_{\vy\in \sss}\lambda\phi(\vy-\vx)\dd\Leb(\vy) \\
  &=&
  \lambda\int_{\vy\in \sss-\vx} \phi(\vy)\dd\Leb(\vy).
\end{eqnarray*}
Since $\int_{\norm{\vy}\le K} \phi(\vy)\to \int_{\RR^d} \phi(\vy)=1$ as $K\to\infty$,
there is a constant $K$ such that
$(T_\kk f)(\vx)\ge (1+\lambda)/2>1$ if $\vx\in \sss$ is at distance at least $K$ from
the boundary of $\sss$.
It follows that if $L$ is large enough, then
$\norm{T_\kk f}_2>\norm{f}_2$, so
$\norm{T_\kk}>1$. From now on we choose $L$ large enough that $\norm{T_\kk}>1$ for $\sss=\sss_1,\sss_2$.

Setting $\rho=r^d$, let $V_\rho=(X/r)\cap \sss$ be the set
of vertices of $G'(X)$ lying in $\sss$. The graph $G_\rho(\kk)$
considered in Theorem~\ref{th_k} has
exactly the distribution of $G'[\sss]$, the subgraph of $G'(X)$ induced by vertices
in $\sss$.
Hence, taking $\sss=\sss_1$ and applying Theorem \ref{th_k},
there is a constant $a=a(\sss_1,\Leb,\kk)>0$ such that,
as $r\to\infty$,
\begin{equation}\label{c1}
 \Pr\bb{ C_1(G'[\sss_1]) \le a\rho } \to 0.
\end{equation}
Taking $\sss=\sss_2$ and applying Theorem \ref{th_k} again, we have
\[
 \Pr\bb{ C_2(G'[\sss_2]) \ge a\rho }  \to 0.
\]

For each $\vv=(v_1,\ldots,v_d)\in \ZZ^d$, let $\sss_\vv=\prod_{i=1}^d [v_iL,v_iL+L)$.
Also, for each edge $e=\{\vv,\vw\}$ of $\ZZ^d$, let $\sss_e=\sss_\vv\cup\sss_\vw$.
We claim that $\Pr\bb{ C_1(G'[\sss_\vv]) \le a\rho }\to 0$ uniformly in $\vv$.
In the Poisson case, this is immediate from \eqref{c1}, since the translation invariance of the model
implies that the relevant probability is independent of $\vv$.
In the general case, a little technical argument is needed: taking $\vv_r$
as any sequence of points of $\ZZ^d$, and defining $V_\rho$,  $\rho=r^d$, by
translating $(X/r)\cap \sss_{\vv_r}$ through $-L\vv_r$, our density assumption
on $X$ implies that $(\sss_1,\bpm,(V_\rho))$ is a vertex space.
Hence, $\Pr\bb{C_1(G'[\sss_{\vv_r}])\le a\rho}\to 0$. As the sequence $\vv_r$
is arbitrary, this is the same as uniform convergence.
Similarly, $\Pr\bb{C_2(G'[\sss_e])\ge a\rho} \to 0$ uniformly in the edges
$e$ of $\ZZ^d$.

For each edge $e=\{\vv,\vw\}$ of the graph $\ZZ^d$, let $\good(e)$ be the event
that $C_1(G'[\sss_\vv]), C_1(G'[\sss_\vw])> a\rho$, while
$C_2(G'[\sss_e])<a\rho$. 
We have shown that $\Pr(\good(e))\to 1$ uniformly in $e$ as $r\to \infty$.
Let $p_0=p_0(d+1)$ be the constant in Lemma~\ref{l_d}, taking $k=d+1$, and choose
$r$ large enough that $\Pr(\good(e))\ge p_0$ for every edge $e$ of $\ZZ^d$.
We shall also assume that $r>D/L$, where $D$ is the constant appearing in
our independence assumption on $X$.

Define a bond percolation measure on $\ZZ^d$ by declaring the edge $e$
to be open if $\good(e)$ holds.
For edges $e$, $f$ at graph distance at least $d+1$,
the sets $\sss_e$ and $\sss_f$ are separated by a Euclidean distance of at least $L>D/r$.
Hence, from our assumptions on $X$ and the independence of edges in the graph,
the graphs $G'[\sss_e]$ and $G'[\sss_f]$ are independent,
and so are the
events $\good(e)$ and $\good(f)$. This observation also holds for sets of edges
at graph distance at least $d+1$, so the bond percolation measure we have
defined is $(d+1)$-independent. Hence, by Lemma~\ref{l_d}, with probability $1$
there is an infinite open path, i.e., an infinite sequence
$\vv_1,\vv_2,\ldots$ such that $\good(\{\vv_i,\vv_{i+1}\})$ holds for each $i\ge 1$.

Let $\CC_i$ be the largest component of $G'[\sss_{\vv_i}]$, chosen arbitrarily 
if there is a tie. As $\good(e)$ holds for $e=\{\vv_i,\vv_{i+1}\}$,
$\CC_i$ and $\CC_{i+1}$ have more than $a\rho$ vertices.
As $G'[\sss_{\vv_j}]$ is a subgraph of $G'[\sss_e]$ for $j=i$, $i+1$,
each of $\CC_i$, $\CC_{i+1}$ is contained entirely within some component
of $G'[\sss_e]$. But as $\good(e)$ holds, $G'[\sss_e]$ has at most one component
with more than $a\rho$ vertices. Hence $\CC_i$ and $\CC_{i+1}$ are connected
in $G'[\sss_e]$, and thus in $G'(X)$. It follows that, with probability 1,
$G'(X)$ and hence $G(X)$ contains an infinite path, completing the proof of Theorem~\ref{th_1}.
\end{proof}

The idea of combining local information (here the events $Y(e)$) to
deduce global information, in particular via comparison with a
product measure, is natural and has been used many times. Here the
events we use are as in Balister, Bollob\'as and
Walters~\cite{BBW}. For an earlier application of related ideas in a
more complicated context, see Pisztora~\cite{P}.

\section{Discussion}

\refT{th_1} is the main part of the related results of
Penrose~\cite{Penrose}, showing that the threshold $\lambda(r)$ for
percolation to occur in $G(X)=G_{r,\lambda}(X)$ approaches 1 as
$r\to\infty$. Note that for $X=\ZZ^d$ or $X$ Poisson, the cases
considered in \cite{Penrose}, the existence of a threshold for each
$r$ follows from monotonicity of the model in $\lambda$ and Kolmogorov's
$0$-$1$ law: constructing $X$ and then $G(X)$ from appropriate
independent random variables, the event that $G(X)$ has an infinite
component is a tail event, and so has probability $0$ or $1$ for any
fixed $\lambda$. Hence there is a (perhaps infinite) $\lambda(r)$ such that this
probability is $0$ for $\lambda<\lambda(r)$ and $1$ for $\lambda>\lambda(r)$.

The condition of \refT{th_1} that $\phi$ be positive in a neighbourhood of
the origin is not essential: it was imposed here for convenience, to avoid
the complication of dealing with reducible kernels in the proof. This
condition is not imposed in \cite{Penrose}. On the other hand,
the stronger conditions on $\phi$ in \cite{Penrose} (or at least some
stronger conditions) are needed for the `easy' part of the result,
that percolation does not occur if $\lambda<1$. This result
is trivial in `nice' cases (see below), but fails under the conditions
of \refT{th_1}, for example if $X=\ZZ^d$ and $\phi$ is large at all
points with integer coordinates
and small elsewhere.

For completeness, we give a short proof of the reverse bound in simple cases;
for proofs under slightly different assumptions see Penrose~\cite{Penrose}
and Meester, Penrose and Sarkar~\cite{MPS}.

\begin{lemma}\label{l_1}
Let $X$ be either $\ZZ^d$ or a Poisson process in $\RR^d$ of intensity $1$,
and let $\phi$ satisfy the conditions of \refT{th_1}. If $X=\ZZ^d$, suppose
in addition that $\phi$ has bounded support.
If $\lambda<1$ is fixed and 
$r$ is large enough, then with probability $1$ every component
of $G(X)$ is finite.
\end{lemma}

\begin{proof}
We start with the Poisson case, working with the rescaled graph $G'(X)$.
Let $U$ be a fixed unit cube in $\RR^d$.
From elementary properties of Poisson processes, the expected number $E_n$ of paths
$(\vx_0,\vx_1,\ldots,\vx_n)$ in $G'(X)$ with $\vx_0\in U$
is given by
\[
 \int_{(\vx_0,\vx_1,\ldots,\vx_n)\in U\times (\RR^d)^n}  r^{d(n+1)} 
    \prod_{i=1}^n\lambda r^{-d}\phi(\vx_i-\vx_{i-1}),
\]
where the integral is with respect to $(d(n+1))$-dimensional Lebesgue
measure, the factor $r^{d(n+1)}$ is due to the density $r^d$ of the (rescaled) Poisson process
on $\RR^d$,
and each factor $\lambda r^{-d}\phi(\vx_i-\vx_{i-1})$ is an edge probability.
As $\int_{\RR^d}\phi(\vx)\dd\Leb(\vx)=1$, we thus have
\[
 E_n = r^{d(n+1)} \Leb(U)\left(\lambda r^{-d}\int_{\RR^d}\phi\right)^n = r^d\lambda^n.
\]
By assumption, $\lambda<1$, so $E_n\to 0$ as $n\to\infty$ with $r$ fixed.
Since a Poisson process has (with probability 1 or by definition) no
accumulation points, every vertex of $G'(X)$ has finite degree. Thus
the probability that $G'(X)$ has an infinite component meeting $U$ is at most
the probability that $G'(X)$ contains a path of length $n$ starting in $U$,
and hence at most $E_n$. So with probability $1$ all components of $G'(X)$ meeting
$U$ are finite. Considering countably many choices for $U$, the result follows.

For the case $X=\ZZ^d$, as $\phi$ is bounded, almost everywhere continuous,
and has bounded support, we have
$r^{-d}\sum_{\vx\in\ZZ^d}\phi(\vx/r)\to \int_{\RR^d}\phi(x)=1$.
(In fact, all we need is that
$\phi$ is directly Riemann integrable, as assumed in \cite{Penrose}.)
Hence, choosing $r$ large enough, we may assume that that
$\lambda r^{-d}\sum_{\vx\in\ZZ^d}\phi(\vx/r)  <1$.
Thus $G(X)$, the unrescaled graph, is a random graph (on $\ZZ^d$) where edges
are independent, and the expected degree of every vertex is at most $c<1$.
Hence the expected number of paths of length $n$ starting at a given vertex
is at most $c^n$, and it follows as above that with probability 1 every component
of $G(X)$ is finite.
\end{proof}

It follows from \refT{th_1}, \refL{l_1} and the remarks on the existence of a threshold
that, under the assumptions of \refL{l_1}, the critical value $\lambda(r)$ for percolation
to occur tends to $1$ as $r\to\infty$. A slightly more general version of this result
is the main result of \cite{Penrose}.

Penrose~\cite{Penrose} also shows that with $\lambda>1$ fixed and $r\to\infty$, 
the probability that the origin (added as an extra point if $X$ is Poisson)
is in an infinite component tends to $\psi(\lambda)$, the survival probability
of a Galton-Watson branching process in which each particle has a Poisson
number of children with mean $\lambda$. The lower bound in this result
also follows from
the results of \cite{kernels}, which relate the size of the giant component
in an inhomogeneous random graph to a branching process. Again, the upper
bound requires stronger conditions.

We close by noting that the results of Penrose considered here are
similar to, but distinct from,
results for `annulus percolation' due
to Balister, Bollob\'as and Walters~\cite{BBW_ann} and Franceschetti,
Booth, Cook, Meester and Bruck~\cite{FBCMB}.  In both cases, the planar
percolation process locally looks like a tree with constant average
degree, and the result is that the average degree at the threshold
approaches $1$, but the methods needed to show this are rather
different in the two cases. (A similar comment applies to the results
of Meester, Penrose and Sarkar~\cite{MPS} and of \cite{BBW_ann}
concerning a related model whose
dimension tends to infinity.) There may well be a common generalization
of these results.


\begin{thebibliography}{}

\bibitem{BBW_ann} P.~Balister, B.~Bollob\'as and M.~Walters,
 Continuum percolation with steps in an annulus,
 {\em Ann. Appl. Probab.} {\bf 14} (2004), 1869--1879.

\bibitem{BBW} P.~Balister, B.~Bollob\'as and M.~Walters,
 Continuum percolation with steps in the square or the disc,
 {\em Random Structures and Algorithms} {\bf 26} (2005), 392--403.

\bibitem{kernels} B.~Bollob\'as, S. Janson and O.~Riordan,
 The phase transition in inhomogeneous random graphs,
 to appear in  {\em Random Structures and Algorithms}.
 Preprint available from \webcite{http://www.arxiv.org/abs/math.PR/0504589}

\bibitem{BK} B.~Bollob{\'a}s and Y.~Kohayakawa,
 A note on long-range percolation,
 in {\em Graph theory, combinatorics, and algorithms, (Kalamazoo, MI, 1992)},
 Wiley-Intersci. Publ. (1995), 97--113.

\bibitem{FBCMB}
 M.~Franceschetti, L.~Booth, M.~Cook, R.~Meester and J.~Bruck,
 Continuum percolation with unreliable and spread-out connections,
 {\em J. Statistical Physics} {\bf 118} (2005), 721--734.

\bibitem{LSS} T.M.~Liggett, R.H.~Schonmann and A.M.~Stacey,
 Domination by product measures,
 {\em Annals of Probability} {\bf 25} (1997), 71--95.

\bibitem{MPS} R.~Meester, M.D.~Penrose and A.~Sarkar,
 The random connection model in high dimensions,
 {\em Statistics and Probability Letters} {\bf 35} (1997), 145--153.

\bibitem{Penrose} M.D.~Penrose,
 On the spread-out limit for bond and continuum percolation,
 {\em Annals of Applied Probability} {\bf 3} (1993), 253--276.

\bibitem{P} A.~Pisztora,
 Surface order large deviations for Ising, Potts and percolation models,
 {\em Probab. Theory Related Fields} {\bf 104} (1996), 427--466.

\end{thebibliography}
\end{document}